\documentclass[11pt]{amsart}
\usepackage[T1]{fontenc}
\usepackage{lmodern}
\usepackage{amsmath,amssymb,amsthm,mathtools}
\usepackage[margin=1.15in]{geometry}
\usepackage[dvipsnames]{xcolor}
\usepackage{hyperref}
\hypersetup{colorlinks=true,linkcolor=NavyBlue,citecolor=NavyBlue,urlcolor=NavyBlue}

\newtheorem{theorem}{Theorem}[section]
\newtheorem{lemma}[theorem]{Lemma}
\newtheorem{proposition}[theorem]{Proposition}
\newtheorem{corollary}[theorem]{Corollary}
\theoremstyle{remark}
\newtheorem{remark}[theorem]{Remark}

\newcommand{\cM}{\mathcal M}
\newcommand{\cH}{\mathcal H}
\newcommand{\cQ}{\mathcal Q}
\newcommand{\bC}{\mathbb C}
\newcommand{\rk}{\operatorname{rank}}
\newcommand{\rel}{\operatorname{rel}}
\newcommand{\Res}{\operatorname{Res}}
\newcommand{\ord}{\operatorname{ord}}

\title[Rank and rel in strata]{Rank and Rel in Strata of Abelian Differentials}
\author{Ruocheng Yang}
\address{Department of Mathematics, University of Toronto, Toronto, ON M5S 2E4, Canada}
\email{ruocheng.yang@mail.utoronto.ca}

\begin{document}

\begin{abstract}
We prove that every invariant subvariety $\cM$ in a genus $g$ stratum $\cH(k_1,\ldots,k_n)$ of Abelian differentials satisfies $\rk(\cM)\le\rel(\cM)+1+\sum_i(k_i-1)$. In particular, every invariant subvariety with rel zero in a principal stratum is a Teichm\"uller curve. Combined with work of Eskin, Filip and Wright, this implies that each principal stratum contains only finitely many nonarithmetic Teichm\"uller curves. If $\rk(\cM)<g$, we also show that $\rk(\cM)+\rel(\cM)\le g-1+\lfloor n/2\rfloor$.
\end{abstract}

\maketitle

\section{Introduction}\label{sec:intro}
A \emph{translation surface} is a pair $(X,\omega)$, where $X$ is a compact connected Riemann surface and $\omega$ is a nonzero holomorphic one form. For $g\ge2$ and positive integers $k_1,\ldots,k_n$ with $\sum_i k_i=2g-2$, the \emph{stratum} $\cH(k_1,\ldots,k_n)$ parametrizes genus $g$ translation surfaces whose zeros have orders $k_1,\ldots,k_n$. The \emph{principal stratum} $\cH(1^{2g-2})$ consists of forms with only simple zeros.

Integrating $\omega$ over a basis of $H_1(X,\Sigma;\mathbb Z)$, where $\Sigma$ is the zero set of $\omega$, gives local \emph{period coordinates}. Eskin, Mirzakhani and Mohammadi proved that orbit closures of the natural $\mathrm{GL}^+(2,\mathbb R)$ action are properly immersed suborbifolds defined locally by real linear equations in these coordinates \cite{EM18, EMM15}. Filip proved that these orbit closures are algebraic \cite{FilipMHS}, and we call them \emph{invariant subvarieties}.

Let $\cM$ be an invariant subvariety containing $x=(X,\omega)$, and identify the tangent space $T_x\cM$ with a subspace of $H^1(X,\Sigma;\bC)$ using period coordinates. Let $p\colon H^1(X,\Sigma;\bC)\to H^1(X;\bC)$ be the natural projection and set $W_x=T_x\cM\cap\ker p$. The \emph{rank} and \emph{rel} of $\cM$ are $\rk(\cM)=\tfrac12\dim_{\bC}p(T_x\cM)$ and $\rel(\cM)=\dim_{\bC}W_x$. They do not depend on $x$, and the rank is an integer because $p(T_x\cM)$ is symplectic \cite{AEM17}.

Previous work has classified invariant subvarieties in various ranges of rank. In this paper, we instead prove general inequalities between rank and rel that hold in every stratum. Our first result gives a lower bound on rel in terms of rank and the orders of the zeros.

\begin{theorem}\label{thm:main}
Let $\cM$ be an invariant subvariety in a stratum $\cH(k_1,\ldots,k_n)$. Then
\[
\rk(\cM)\le\rel(\cM)+1+\sum_{i=1}^n(k_i-1).
\]
\end{theorem}

\begin{corollary}\label{cor:relzero}
An invariant subvariety with rel zero in a stratum $\cH(k_1,\ldots,k_n)$ has rank at most $1+\sum_i(k_i-1)$. In particular, every invariant subvariety with rel zero in a principal stratum is a Teichm\"uller curve.
\end{corollary}

Together with a theorem of Eskin, Filip and Wright, Theorem~\ref{thm:main} gives the following finiteness result.

\begin{corollary}\label{cor:finiteTC}
For every $g\ge2$, the principal stratum of genus $g$ contains only finitely many nonarithmetic Teichm\"uller curves.
\end{corollary}

Mirzakhani and Wright classified invariant subvarieties of full rank $g$, showing that they are connected components of strata or hyperelliptic loci \cite[Theorem~1.1]{MW18}. Our second result complements Theorem~\ref{thm:main} by giving an upper bound on rel for invariant subvarieties of rank less than $g$.

\begin{theorem}\label{thm:upper}
Let $\cM$ be an invariant subvariety in a stratum $\cH(k_1,\ldots,k_n)$ of genus $g$. If $\cM$ does not have full rank, then
\[
\rk(\cM)+\rel(\cM)\le g-1+\lfloor n/2\rfloor.
\]
\end{theorem}

Since $n\le2g-2$, Theorem~\ref{thm:upper} implies that $\rk(\cM)+\rel(\cM)\le2g-2$, with equality only possible in principal strata. The principal stratum itself has rank $g$ and rel $2g-3$, so the assumption on the rank cannot be dropped. Theorem~\ref{thm:upper} also gives the following bound for invariant subvarieties of maximal rel.

\begin{corollary}\label{cor:fullrel}
Let $\cM$ be an invariant subvariety of rel $n-1$ in a stratum $\cH(k_1,\ldots,k_n)$ of genus $g$. If $\cM$ does not have full rank, then $\rk(\cM)\le g-\lceil n/2\rceil$. In particular, every proper invariant subvariety of maximal rel in a principal stratum has rank one; if $g\ge3$, every such subvariety is a locus of torus covers.
\end{corollary}

Calsamiglia, Deroin and Francaviglia classified closures of leaves of the absolute period foliation of $\Omega\cM_g$ for $g\ge3$ \cite{CDF23}, and Winsor gave a new proof of the resulting classification of saturated orbit closures \cite{Win24}. In particular, every proper invariant subvariety of maximal rel in the principal stratum of genus $g\ge3$ is a locus of torus covers. Corollary~\ref{cor:fullrel} gives a different proof of this torus cover conclusion, as well as the rank one conclusion in genus two. Winsor has since extended the description of leaf closures to connected strata, away from a proper algebraic subvariety \cite{Win26}. See \cite{CW24, Win21} for the dynamics of the real Rel foliation and the absolute period foliation.

McMullen classified invariant subvarieties in genus two \cite{McM07}, and those of rank at least two in genus three were classified in \cite{ANW16, AN16, AN20}. In arbitrary genus, Apisa and Wright extended the full rank classification to invariant subvarieties of rank at least $g/2+1$ \cite{AW23}, while Apisa classified those of rank at least two in hyperelliptic components \cite{Api18}. In the principal stratum, M\"oller showed that every invariant subvariety other than the stratum and its hyperelliptic locus parametrizes surfaces whose Jacobians have extra endomorphisms \cite{MollerLinear}. Eskin, Filip and Wright showed that each stratum contains only finitely many invariant subvarieties of rank at least two \cite{EFW18}.

\medskip
\noindent\textbf{Acknowledgments.} The author is grateful to Paul Apisa for carefully reading an earlier version of the manuscript and for many valuable suggestions, to Alex Wright for helpful discussions and suggestions concerning the proof, and to Kai Fu for insightful discussions and comments on an earlier draft. The author also thanks David Aulicino and Curt McMullen for their comments on the manuscript.

\medskip
\noindent\textbf{AI disclosure.} The author used ChatGPT for language editing, expository suggestions and checking bibliographic details. The author reviewed the suggestions and takes full responsibility for the content of this paper.

\section{Proof of Theorem~\ref{thm:main}}
Fix $x=(X,\omega)\in\cM$ and work locally on a smooth branch of $\cM$ through $x$. As in the definition of period coordinates, we use the Gauss--Manin connection to identify the relative cohomology of nearby surfaces with $H^1(X,\Sigma;\bC)$, compatibly with $p$. Label the zeros $\Sigma=\{z_1,\ldots,z_n\}$ so that $z_i$ has order $k_i$, and recall the exact sequence
\[
0\longrightarrow\bC^\Sigma/\bC\mathbf1
\longrightarrow H^1(X,\Sigma;\bC)
\xrightarrow{p}H^1(X;\bC)\longrightarrow0,
\]
where $\bC^\Sigma$ is the space of functions on $\Sigma$ and $\bC\mathbf1$ is the subspace of constant functions. We normalize the injection so that the class of $c\colon\Sigma\to\bC$ vanishes on absolute cycles and takes the value $c(z_i)-c(z_1)$ on a path from $z_1$ to $z_i$.

For a holomorphic one form $\alpha$, let $\widetilde\alpha\in H^1(X,\Sigma;\bC)$ be its relative class, so that $\widetilde\alpha(\gamma)=\int_\gamma\alpha$ for $\gamma\in H_1(X,\Sigma;\mathbb Z)$. We identify holomorphic one forms with their classes in $H^{1,0}(X)$, so that $p(\widetilde\alpha)=\alpha$, and write $p(T_x\cM)^{1,0}=p(T_x\cM)\cap H^{1,0}(X)$.

Write $r=\rk(\cM)$ and $s=\rel(\cM)$, so that $\dim\cM=2r+s$. Because $p(T_x\cM)$ is defined over $\mathbb R$ in period coordinates and compatible with the Hodge decomposition \cite{FilipKZ}, we have $\dim p(T_x\cM)^{1,0}=r$. Complex scaling shows that $\omega\in p(T_x\cM)^{1,0}$.

We use the Hodge filtration $F^\bullet$ and the filtrations it induces on subspaces and quotients. In particular, $F^1H^1(X;\bC)=H^{1,0}(X)$, and $F^1H^1(X,\Sigma;\bC)$ consists of the relative classes of holomorphic one forms. We write $\mathcal O_X$, $K_X$ and $T_X$ for the structure sheaf, canonical bundle and tangent bundle of $X$, and $h^0(X,L)=\dim H^0(X,L)$ for a line bundle $L$ on $X$.

The following lemma is a consequence of Filip's splitting theorem \cite[Theorem~4.2]{FilipMHS}. We include the deduction for completeness.

\begin{lemma}\label{lem:Filip}
For every $\alpha\in p(T_x\cM)^{1,0}$, one has $\widetilde\alpha\in T_x\cM$.
\end{lemma}

\begin{proof}
Consider the exact sequence of real mixed Hodge structures, written after complexification,
\[
0\longrightarrow \ker p/W_x \longrightarrow p^{-1}(p(T_x\cM))/W_x \xrightarrow{p} p(T_x\cM) \longrightarrow 0,
\]
where $p^{-1}$ denotes the preimage. Since $\ker p/W_x$ has Hodge type $(0,0)$, its intersection with $F^1(p^{-1}(p(T_x\cM))/W_x)$ is zero, so $p$ is injective on $F^1(p^{-1}(p(T_x\cM))/W_x)$.

The projection $p$ induces an isomorphism $T_x\cM/W_x\to p(T_x\cM)$, whose inverse we denote by
\[
\sigma_x\colon p(T_x\cM)\longrightarrow T_x\cM/W_x \subset p^{-1}(p(T_x\cM))/W_x.
\]
Filip's splitting theorem \cite[Theorem~4.2]{FilipMHS} implies that $\sigma_x$ preserves the Hodge filtration, so $\sigma_x(\alpha)\in F^1(p^{-1}(p(T_x\cM))/W_x)$ for every $\alpha\in p(T_x\cM)^{1,0}=F^1(p(T_x\cM))$.

On the other hand, the relative class $\widetilde\alpha$ of such an $\alpha$ lies both in $F^1H^1(X,\Sigma;\bC)$ and in $p^{-1}(p(T_x\cM))$, so $\widetilde\alpha+W_x\in F^1(p^{-1}(p(T_x\cM))/W_x)$. Since $\widetilde\alpha+W_x$ and $\sigma_x(\alpha)$ both project to $\alpha$, injectivity gives
\[
\widetilde\alpha+W_x=\sigma_x(\alpha)\in T_x\cM/W_x,
\]
and hence $\widetilde\alpha\in T_x\cM$ because $W_x\subset T_x\cM$.
\end{proof}

Let $\pi\colon\cH(k_1,\ldots,k_n)\to\cM_g$ be the forgetful map $(X,\omega)\mapsto X$. We now describe $\ker d\pi_x$. Working in an orbifold chart, we regard the stratum near $x$ as a complex submanifold of the Hodge bundle $\Omega\cM_g$, so that $\pi$ is the restriction of the bundle projection $\Omega\cM_g\to\cM_g$. The tangent space to the fiber of $\Omega\cM_g$ over $X$ is naturally $H^0(X,K_X)$, and we call its elements \emph{vertical tangent vectors}. Thus $\ker d\pi_x$ consists of the vertical tangent vectors that are tangent to the stratum. Let $D=\sum_i(k_i-1)z_i$, so that $H^0(X,K_X(-D))$ consists of the holomorphic one forms $\beta$ with $\ord_{z_i}\beta\ge k_i-1$ for all $i$, and set
\[
V_x=p(T_x\cM)^{1,0}\cap H^0(X,K_X(-D)).
\]

\begin{lemma}\label{lem:kerdpi}
One has
\[
\ker d\pi_x=\{\widetilde\beta:\beta\in H^0(X,K_X(-D))\}
\quad\text{and}\quad
\ker d(\pi|_\cM)_x=\{\widetilde\beta:\beta\in V_x\}.
\]
In particular, $\dim\ker d(\pi|_\cM)_x=\dim V_x$.
\end{lemma}

\begin{proof}
We view a vertical tangent vector $\beta\in H^0(X,K_X)$ as the first-order deformation $\omega+\varepsilon\beta$ of $\omega$ over $\bC[\varepsilon]/(\varepsilon^2)$, with $X$ fixed. By \cite[Theorems~2.1 and~2.3]{MollerLinear} (see also \cite[Theorem~2.1]{BCGGM19}), tangent vectors to the stratum at $x$ are exactly the first-order deformations of $(X,\omega,z_1,\ldots,z_n)$ in which each $z_i$ deforms to a zero of order $k_i$. So $\beta$ is tangent to the stratum if and only if each $z_i$ deforms to a zero of order $k_i$ of $\omega+\varepsilon\beta$. Away from $\Sigma$ the coefficient of $\omega+\varepsilon\beta$ remains invertible, so no new zeros appear there.

To check the condition at $z_i$, pick a local coordinate $z$ centered at $z_i$ with $\omega=z^{k_i}f(z)\,dz$ and $f(0)\ne0$, and write $\beta=b(z)\,dz$. A first-order deformation of $z_i$ has velocity $v\in\bC$ and is locally defined by $z-\varepsilon v=0$. So $z_i$ deforms to a zero of order $k_i$ of $\omega+\varepsilon\beta$ exactly when $z^{k_i}f+\varepsilon b$ is $(z-\varepsilon v)^{k_i}$ times a unit. Setting $\varepsilon=0$ shows that this unit reduces to $f$, so it has the form $f+\varepsilon h$ for a holomorphic function $h$ near $0$, and conversely every such function is a unit because $f(0)\ne0$. Since $\varepsilon^2=0$, we also have $(z-\varepsilon v)^{k_i}=z^{k_i}-\varepsilon k_iv\,z^{k_i-1}$. Thus the condition is that
\[
z^{k_i}f+\varepsilon b=\bigl(z^{k_i}-\varepsilon k_iv\,z^{k_i-1}\bigr)(f+\varepsilon h)
\]
for some $v\in\bC$ and some holomorphic function $h$ near $0$. Comparing the coefficients of $\varepsilon$, this means that
\[
b=-k_iv\,z^{k_i-1}f+z^{k_i}h.
\]
The right side vanishes to order at least $k_i-1$, so $\ord_{z_i}\beta\ge k_i-1$ is necessary. Conversely, if $\ord_{z_i}\beta\ge k_i-1$, write $b=c\,z^{k_i-1}+O(z^{k_i})$. Then $b+k_iv\,z^{k_i-1}f$ vanishes to order at least $k_i$ exactly when $v=-c/(k_if(0))$, and for this $v$ the function $h=(b+k_iv\,z^{k_i-1}f)/z^{k_i}$ is holomorphic. Applying this at every $z_i$ shows that the vertical tangent vectors to the stratum are exactly the elements of $H^0(X,K_X(-D))$.

It remains to identify these vectors in period coordinates. For an absolute cycle $\gamma$, the first-order variation of its period is $\int_\gamma\beta$. For a path $\gamma_i$ from $z_1$ to $z_i$, we must also account for the motion of its endpoints. Let $v_j$ be the velocity of $z_j$ found above, and write $\omega=z^{k_j}f_j(z)\,dz$ in the coordinate centered at $z_j$. Extending $\gamma_i$ by short paths to the deformed zeros, we find that modulo $\varepsilon^2$ the period of $\omega+\varepsilon\beta$ along the extended path is
\[
\int_{\gamma_i}\omega+\varepsilon\int_{\gamma_i}\beta+\int_0^{\varepsilon v_i}z^{k_i}f_i\,dz-\int_0^{\varepsilon v_1}z^{k_1}f_1\,dz,
\]
since the integrals of $\varepsilon\beta$ along the short paths are multiples of $\varepsilon^2$. The last two integrals are multiples of $\varepsilon^{k_i+1}$ and $\varepsilon^{k_1+1}$, so they vanish modulo $\varepsilon^2$. Hence the first-order variation of the relative period is $\int_{\gamma_i}\beta$. Since absolute cycles and the paths $\gamma_i$ generate $H_1(X,\Sigma;\mathbb Z)$, the corresponding tangent vector is $\widetilde\beta$, proving the first equality.

For the second equality, note that $\ker d(\pi|_\cM)_x=T_x\cM\cap\ker d\pi_x$, which by the first equality consists of the classes $\widetilde\beta$ with $\beta\in H^0(X,K_X(-D))$ and $\widetilde\beta\in T_x\cM$. For such $\beta$, the condition $\widetilde\beta\in T_x\cM$ is equivalent to $\beta\in p(T_x\cM)^{1,0}$: if $\widetilde\beta\in T_x\cM$, then $\beta=p(\widetilde\beta)$ lies in $p(T_x\cM)^{1,0}$, and the converse is Lemma~\ref{lem:Filip}. Thus $\ker d(\pi|_\cM)_x=\{\widetilde\beta:\beta\in V_x\}$, and the last assertion follows because $\beta\mapsto\widetilde\beta$ is injective.
\end{proof}

This identification of the vertical tangent space can also be read off from \cite[Theorem~2.3]{MollerLinear}, which describes the tangent space of the stratum at $x$ as the first hypercohomology of the complex $T_X\to K_X(-D)$ given by $v\mapsto d(\omega(v))$. In this description $d\pi_x$ is the natural map to $H^1(X,T_X)$, and its kernel is identified with $H^0(X,K_X(-D))$ because $H^0(X,T_X)=0$ for $g\ge2$.

By Lemma~\ref{lem:kerdpi}, $\dim\ker d(\pi|_\cM)_x=\dim V_x$, which may vary with $x$ (cf.~\cite[Section~2.2]{MollerLinear}). The next lemma requires this dimension to be constant near $x$. A nonvanishing minor of $d(\pi|_\cM)$ remains nonvanishing nearby, so the maximal rank locus of $d(\pi|_\cM)$ is a nonempty open subset of the branch, on which $\dim\ker d(\pi|_\cM)$ is constant. Since $r$, $s$ and the $k_i$ do not depend on $x$, we may assume in the proof of Theorem~\ref{thm:main} that $x$ lies in this locus.

\begin{lemma}\label{lem:fixedfiber}
Let $\alpha\in V_x$. Then there is a holomorphic family $x_t=(X,\omega_t)$ in $\cM$, defined for small $t\in\bC$, with $x_0=x$ and tangent vector $\widetilde\alpha$ at $t=0$. Moreover, for all small $t$, the relative class $\widetilde\alpha_t$ of $\alpha$ with respect to the zeros of $\omega_t$ lies in $T_{x_t}\cM$.
\end{lemma}

\begin{proof}
Because the rank of $d(\pi|_\cM)$ is constant near $x$, the holomorphic constant rank theorem \cite[Theorem~2.3]{Lee24} shows that a neighborhood of $x$ in the fiber of $\pi|_\cM$ over $X$ is a complex submanifold $F$ of $\cM$ with $T_xF=\ker d(\pi|_\cM)_x$. This tangent space contains $\widetilde\alpha$ by Lemma~\ref{lem:kerdpi}, so there is a holomorphic curve $x_t$ in $F$ with $x_0=x$ and tangent vector $\widetilde\alpha$ at $t=0$. It has the form $x_t=(X,\omega_t)$ because $F$ lies over $X$.

Along this curve, $T_{x_t}\cM=T_x\cM$ under the Gauss--Manin identification, as $\cM$ is locally linear in period coordinates. The absolute cohomology and its Hodge decomposition do not change either, as $X$ is fixed. So $\alpha\in p(T_{x_t}\cM)^{1,0}$, and Lemma~\ref{lem:Filip} applied at $x_t$ gives $\widetilde\alpha_t\in T_{x_t}\cM$. Note that this does not require $\alpha\in V_{x_t}$, which may fail for $t\ne0$ since the zeros of $\omega_t$ move.
\end{proof}

Let $\alpha\in V_x$. The meromorphic one form $\alpha^2/\omega$ is holomorphic away from $\Sigma$, and at each $z_i$
\[
\ord_{z_i}\frac{\alpha^2}{\omega}\ge2(k_i-1)-k_i=k_i-2.
\]
So $\alpha^2/\omega$ has at most simple poles, all at simple zeros of $\omega$. In particular, its residue vanishes at each $z_i$ with $k_i\ge2$. We define the \emph{quadratic residue map} $Q_\omega\colon V_x\to\bC^\Sigma/\bC\mathbf1$ by
\[
Q_\omega(\alpha)=\left[\Res_{z_1}\frac{\alpha^2}{\omega},\ldots,\Res_{z_n}\frac{\alpha^2}{\omega}\right],
\]
where the brackets denote the class modulo constant functions.

\begin{proposition}\label{prop:inW}
For every $\alpha\in V_x$, one has $Q_\omega(\alpha)\in W_x$.
\end{proposition}

\begin{proof}
Fix $\alpha\in V_x$ and let $x_t=(X,\omega_t)$ be the family from Lemma~\ref{lem:fixedfiber}. Since $X$ is fixed, the tangent vector of this family at $t=0$ is the vertical tangent vector $\left.\frac{d}{dt}\right|_{t=0}\omega_t$, and the proof of Lemma~\ref{lem:kerdpi} shows that in period coordinates it is the relative class of this form. By Lemma~\ref{lem:fixedfiber}, this tangent vector is $\widetilde\alpha$, so the injectivity of $\beta\mapsto\widetilde\beta$ gives $\left.\frac{d}{dt}\right|_{t=0}\omega_t=\alpha$, that is, $\omega_t=\omega+t\alpha+O(t^2)$.

The zeros of $\omega_t$ move holomorphically. Indeed, choose disjoint coordinate discs around the zeros $z_i$, and write $\omega_t=F_i(t,z)\,dz$ in the disc around $z_i$. By Rouch\'e's theorem, for all small $t$ the form $\omega_t$ has $k_i$ zeros in this disc, counted with multiplicity. Since $x_t$ lies in the stratum, $\omega_t$ has exactly $n$ distinct zeros, so each disc contains exactly one zero of $\omega_t$, and it has order $k_i$. On the other hand, $\partial_z^{k_i-1}F_i(0,z)$ has a simple zero at $z_i$, so the holomorphic implicit function theorem gives a holomorphic function $z_i(t)$ with $z_i(0)=z_i$ such that, for small $t$, $z_i(t)$ is the only solution near $z_i$ of $\partial_z^{k_i-1}F_i(t,z)=0$. The zero of $\omega_t$ in the disc around $z_i$ also satisfies this equation, so it is $z_i(t)$. Thus $z_i(t)$ is a zero of $\omega_t$ of order $k_i$.

By Lemma~\ref{lem:fixedfiber}, $\widetilde\alpha_t\in T_{x_t}\cM=T_x\cM$ for all small $t$. Moreover $p(\widetilde\alpha_t)=\alpha$ for all $t$, so
\[
\left.\frac{d}{dt}\right|_{t=0}\widetilde\alpha_t\in T_x\cM\cap\ker p=W_x.
\]

To compute this derivative, fix $2\le i\le n$, and let $\gamma_i(t)$ be the path from $z_1(t)$ to $z_i(t)$ obtained by transporting a path $\gamma_i(0)$ from $z_1$ to $z_i$. Since $\alpha$ is closed, Stokes' theorem gives
\[
\int_{\gamma_i(t)}\alpha-\int_{\gamma_i(0)}\alpha=\int_{z_i}^{z_i(t)}\alpha-\int_{z_1}^{z_1(t)}\alpha,
\]
where the last two integrals are along the paths traced out by the endpoints. Write $\alpha=a_j(z)\,dz$ in a local coordinate $z$ near each $z_j$. Differentiating gives
\[
\left.\frac{d}{dt}\right|_{t=0}\int_{\gamma_i(t)}\alpha=a_i(z_i)z_i'(0)-a_1(z_1)z_1'(0).
\]
By definition of the transport, $\widetilde\alpha_t(\gamma_i(0))=\int_{\gamma_i(t)}\alpha$, and the values of $\widetilde\alpha_t$ on absolute cycles do not depend on $t$. Since absolute cycles and the paths $\gamma_i(0)$ generate relative homology, our convention for $\bC^\Sigma/\bC\mathbf1$ gives
\[
\left.\frac{d}{dt}\right|_{t=0}\widetilde\alpha_t=\left[a_1(z_1)z_1'(0), \ldots, a_n(z_n)z_n'(0)\right].
\]

If $k_i\ge2$, the $i$-th entry of this vector and $\Res_{z_i}(\alpha^2/\omega)$ both vanish, since $a_i(z_i)=0$ for $\alpha\in H^0(X,K_X(-D))$ and $\ord_{z_i}(\alpha^2/\omega)\ge k_i-2\ge0$. If $k_i=1$, write $\omega=(z-z_i)f_i(z)\,dz$ with $f_i(z_i)\ne0$ in the same coordinate. Evaluating the coefficient of $\omega_t$ at its zero $z_i(t)$ gives
\[
(z_i(t)-z_i)f_i(z_i(t))+ta_i(z_i(t))+O(t^2)=0,
\]
and differentiating at $t=0$ gives $f_i(z_i)z_i'(0)+a_i(z_i)=0$. Thus
\[
a_i(z_i)z_i'(0)=-\frac{a_i(z_i)^2}{f_i(z_i)}=-\Res_{z_i}\frac{\alpha^2}{\omega}.
\]
So
\[
\left.\frac{d}{dt}\right|_{t=0}\widetilde\alpha_t=-Q_\omega(\alpha),
\]
and since this derivative is in $W_x$, so is $Q_\omega(\alpha)$.
\end{proof}

\begin{proposition}\label{prop:zeroset}
The map $Q_\omega\colon V_x\to\bC^\Sigma/\bC\mathbf1$ is homogeneous quadratic. Its zero set is the linear subspace of $V_x$ consisting of the forms that vanish at every simple zero of $\omega$, and this subspace has dimension at most $h^0(X,\mathcal O_X(D))$.
\end{proposition}

\begin{proof}
The map is homogeneous quadratic by definition. If $Q_\omega(\alpha)=0$, then all residues of $\alpha^2/\omega$ are equal, hence zero by the residue theorem. At a simple zero $z_i$ of $\omega$, the form $\alpha^2/\omega$ has order $2\ord_{z_i}\alpha-1$, so its residue there vanishes if and only if $\alpha$ vanishes at $z_i$. Conversely, if $\alpha\in V_x$ vanishes at every simple zero of $\omega$, then $\alpha^2/\omega$ is holomorphic there, and also at the zeros with $k_i\ge2$, where its order is at least $k_i-2\ge0$, so $Q_\omega(\alpha)=0$. Thus the zero set of $Q_\omega$ is the linear subspace of $V_x$ consisting of the forms that vanish at every simple zero of $\omega$.

For $\alpha$ in this subspace, the meromorphic function $\alpha/\omega$ is holomorphic away from $\Sigma$ and has order $\ord_{z_i}\alpha-k_i$ at each $z_i$. This order is at least $0$ when $k_i=1$, since $\alpha$ vanishes there, and at least $-1\ge-(k_i-1)$ when $k_i\ge2$, since $\ord_{z_i}\alpha\ge k_i-1$. Hence $\alpha/\omega\in H^0(X,\mathcal O_X(D))$, and since the map $\alpha\mapsto\alpha/\omega$ is linear and injective, the dimension bound follows.
\end{proof}

\begin{proof}[Proof of Theorem~\ref{thm:main}]
Since $V_x$ is the intersection of $p(T_x\cM)^{1,0}$ and $H^0(X,K_X(-D))$ inside $H^{1,0}(X)$, the dimension formula and Riemann--Roch give
\[
\begin{aligned}
\dim V_x
&\ge r+h^0(X,K_X(-D))-g\\
&=r-1-\deg D+h^0(X,\mathcal O_X(D)).
\end{aligned}
\]

By Proposition~\ref{prop:inW}, the quadratic map $Q_\omega$ takes values in the $s$-dimensional space $W_x$, so after choosing a basis of $W_x$, its zero set is the common zero set of $s$ homogeneous quadratic polynomials on $V_x$. This zero set contains $0$, so it has dimension at least $\dim V_x-s$ \cite[Exercise~I.1.9]{Hartshorne}. On the other hand, it has dimension at most $h^0(X,\mathcal O_X(D))$ by Proposition~\ref{prop:zeroset}. Consequently,
\[
r-1-\deg D+h^0(X,\mathcal O_X(D)) \le\dim V_x \le s+h^0(X,\mathcal O_X(D)),
\]
so $r\le s+1+\deg D$, which is Theorem~\ref{thm:main} since $\deg D=\sum_i(k_i-1)$.
\end{proof}

\begin{proof}[Proof of Corollary~\ref{cor:relzero}]
The first assertion follows from Theorem~\ref{thm:main}. In a principal stratum $D=0$, so if $s=0$, then $r\le1$ by Theorem~\ref{thm:main}, and $r=1$ since $\omega\in p(T_x\cM)^{1,0}$. Thus $\dim\cM=2$, and $\cM$ is a Teichm\"uller curve \cite{Wright15}.
\end{proof}

\begin{proof}[Proof of Corollary~\ref{cor:finiteTC}]
Suppose that the principal stratum contains infinitely many nonarithmetic Teichm\"uller curves. By \cite[Theorem~1.5]{EFW18} and \cite{Wri14}, one of the corresponding closed orbits is contained with complex codimension two in an invariant subvariety $\mathcal N$ of rank two. Since the orbit has complex dimension two, $\dim\mathcal N=4$ and hence $\rel(\mathcal N)=0$, contradicting Theorem~\ref{thm:main}.
\end{proof}

A different argument, using Abel's theorem, gives a stronger bound on the fibers of the forgetful map and a refinement of Theorem~\ref{thm:main}.

\begin{proposition}\label{prop:fiberdim}
Every fiber of $\pi|_\cM$ has dimension at most $s+1$.
\end{proposition}

\begin{proof}
Fix $x=(X,\omega)\in\cM$, and let $F$ be the germ at $x$ of the fiber of $\pi|_\cM$ on a smooth local branch of $\cM$. For $x'=(X,\omega')\in F$, write $z_i'$ for the zero of $\omega'$ near $z_i$, as in the proof of Proposition~\ref{prop:inW}. Since $X$ is fixed and $p(T\cM)$ is flat, $\omega\in p(T_{x'}\cM)^{1,0}$. By Lemma~\ref{lem:Filip}, the relative class of $\omega$ at $x'$ belongs to $T_{x'}\cM=T_x\cM$. Subtracting its relative class at $x$ gives a holomorphic map
\[
F\longrightarrow W_x,\qquad
x'\longmapsto
\left[\left(\int_{z_i}^{z_i'}\omega\right)_{i=1}^n\right],
\]
where the integrals are taken in fixed coordinate discs around the $z_i$. We choose each disc so small that $z_i$ is the only zero of $z\mapsto\int_{z_i}^z\omega$ in it, which is possible because this zero is isolated.

We claim that its zero fiber is precisely the local scaling orbit of $\omega$. The divisors $\sum_ik_iz_i'$ and $\sum_ik_iz_i$ are both canonical, hence linearly equivalent, so by Abel's theorem the functional $\eta\mapsto\sum_ik_i\int_{z_i}^{z_i'}\eta$ on $H^0(X,K_X)$ lies in the period lattice of the Jacobian. This functional tends to zero as $x'$ tends to $x$, and the period lattice is discrete, so it vanishes for $x'$ sufficiently close to $x$. In particular, $\sum_ik_i\int_{z_i}^{z_i'}\omega=0$. If the map vanishes at $x'$, all the integrals are equal, and hence vanish since $\sum_ik_i=2g-2\ne0$. By the choice of the discs, this implies $z_i'=z_i$ for all $i$, so $\omega'/\omega$ is a holomorphic function on $X$, hence constant.

The zero fiber therefore has dimension one. Since $\dim W_x=s$, the local dimension inequality gives $\dim_xF\le s+1$. As $x$ was arbitrary, the assertion follows.
\end{proof}

At a point $x$ where $d(\pi|_\cM)$ has maximal rank, Proposition~\ref{prop:fiberdim} and Lemma~\ref{lem:kerdpi} give $\dim V_x\le s+1$. Since $K_X(-D)\simeq\mathcal O_X(z_1+\cdots+z_n)$, combining this with the lower bound $\dim V_x\ge r+h^0(X,K_X(-D))-g$ from the proof of Theorem~\ref{thm:main} yields
\[
r\le s+g+1-h^0(X,\mathcal O_X(z_1+\cdots+z_n)).
\]

The two arguments are related. At a point $x$ where $d(\pi|_\cM)$ has maximal rank, let $\alpha\in V_x$ and let $x_t=(X,\omega_t)$ be the family from Lemma~\ref{lem:fixedfiber}, so that $\omega_t=\omega+t\alpha+O(t^2)$. At a simple zero $z_i$, with notation as in the proof of Proposition~\ref{prop:inW}, we have $\int_{z_i}^{z_i(t)}\omega=\frac12f_i(z_i)\bigl(z_i(t)-z_i\bigr)^2+O(t^3)$ and $z_i'(0)=-a_i(z_i)/f_i(z_i)$, while at a zero of order $k_i\ge2$ the integral is $O(t^3)$. Hence the map constructed above has expansion
\[
\left[\left(\int_{z_i}^{z_i(t)}\omega\right)_{i=1}^n\right]=\frac{t^2}{2}Q_\omega(\alpha)+O(t^3),
\]
so $Q_\omega$ is the quadratic form associated with the Hessian of this map at $x$.

\section{Proof of Theorem~\ref{thm:upper}}
Fix an arbitrary point $x=(X,\omega)\in\cM$ and work locally on a smooth branch through $x$, keeping the notation of the previous section. We prove Theorem~\ref{thm:upper} by bounding separately the kernel and image of $d(\pi|_\cM)_x$.

We use de Rham models for the following two cohomology groups. By the de Rham theorem, $H^1(X,\Sigma;\bC)$ is the space of closed smooth one forms on $X$ vanishing near $\Sigma$, modulo $d\psi$ for smooth functions $\psi$ vanishing near $\Sigma$, where such a form $\theta$ gives the class $\gamma\mapsto\int_\gamma\theta$ on $H_1(X,\Sigma;\mathbb Z)$. Similarly, $H^1(X\setminus\Sigma;\bC)$ is the space of closed smooth one forms on $X\setminus\Sigma$, with no condition near $\Sigma$, modulo exact forms. In particular, meromorphic one forms on $X$ with poles in $\Sigma$, which are holomorphic and hence closed on $X\setminus\Sigma$, define classes in $H^1(X\setminus\Sigma;\bC)$. For a closed form $\nu$, we write $[\nu]$ for its de Rham class, and we use the same notation for the Dolbeault class in $H^1(X,\mathcal O_X)$ of a $(0,1)$-form, that is, its class modulo $\bar\partial g$ for smooth functions $g$.

For $\xi\in H^1(X,\Sigma;\bC)$ represented by $\theta$ and $[\nu]\in H^1(X\setminus\Sigma;\bC)$, set
\[
\langle\xi,[\nu]\rangle=\frac{1}{2\pi i}\int_{X\setminus\Sigma}\nu\wedge\theta.
\]
The integral is well defined because $\theta$ has compact support in $X\setminus\Sigma$, even though $\nu$ may be unbounded near $\Sigma$. It does not depend on the representatives: if $\psi$ is a smooth function vanishing near $\Sigma$, then $\nu\wedge d\psi=-d(\psi\nu)$, and if $\varphi$ is a smooth function on $X\setminus\Sigma$, then $d\varphi\wedge\theta=d(\varphi\theta)$. In both cases the form in parentheses has compact support in $X\setminus\Sigma$, so its differential integrates to zero by Stokes' theorem. Under the natural identification of $H^1(X,\Sigma;\bC)$ with $H^1_c(X\setminus\Sigma;\bC)$, this is the Poincar\'e duality pairing up to normalization.

By de Rham's theorem, the cup product pairing on $H^1(X;\bC)$, obtained by evaluating the cup product on the fundamental class of $X$, is given by $([a],[b])\mapsto\int_Xa\wedge b$ for closed one forms $a$ and $b$.

\begin{lemma}\label{lem:pairing}
Let $\xi\in H^1(X,\Sigma;\bC)$ be represented by $\theta$ as above. If $\nu$ extends to a closed one form on $X$, then $\langle\xi,[\nu]\rangle$ is $(2\pi i)^{-1}$ times the cup product pairing of $[\nu]$ with $p(\xi)$, namely
\[
\langle\xi,[\nu]\rangle=\frac{1}{2\pi i}\int_X\nu\wedge\theta.
\]
If $\xi\in\ker p$ corresponds to $[c]\in\bC^\Sigma/\bC\mathbf1$ and $\nu$ is meromorphic with poles in $\Sigma$, then
\[
\langle\xi,[\nu]\rangle=\sum_{i=1}^n c(z_i)\Res_{z_i}\nu.
\]
\end{lemma}

\begin{proof}
For the first assertion, $\theta$ vanishes near $\Sigma$, so the integral defining the pairing may be taken over $X$. Moreover, $\theta$ is a closed form on $X$ representing $p(\xi)$ in $H^1(X;\bC)$, since $p$ is induced by the inclusion of absolute cycles into relative cycles. As $\nu$ represents $[\nu]$, the integral of $\nu\wedge\theta$ over $X$ is the cup product pairing of $[\nu]$ with $p(\xi)$.

For the second, choose a smooth function $\varphi$ equal to $c(z_i)$ near each $z_i$, so that $d\varphi$ represents $\xi$. Let $U_i$ be pairwise disjoint small discs around $z_i$ on which $\varphi$ is constant, with $\partial U_i$ oriented as the boundary of $U_i$, and set $X^\circ=X\setminus\bigcup_i\operatorname{int}(U_i)$. Since $d\varphi$ vanishes on the $U_i$ and $\nu\wedge d\varphi=-d(\varphi\nu)$ on $X^\circ$, Stokes' theorem gives
\[
\begin{aligned}
\langle\xi,[\nu]\rangle
&=\frac{1}{2\pi i}\int_{X^\circ}\nu\wedge d\varphi
=-\frac{1}{2\pi i}\int_{\partial X^\circ}\varphi\nu\\
&=\frac{1}{2\pi i}\sum_{i=1}^n c(z_i)\int_{\partial U_i}\nu
=\sum_{i=1}^n c(z_i)\Res_{z_i}\nu.
\end{aligned}
\]
The last expression is independent of the representative $c$, since the residues of $\nu$ sum to zero.
\end{proof}

Let $p(T_x\cM)^\perp$ be the orthogonal complement of $p(T_x\cM)$ in $H^1(X;\bC)$ with respect to the cup product pairing. Since $p(T_x\cM)$ is symplectic \cite{AEM17}, we have $H^1(X;\bC)=p(T_x\cM)\oplus p(T_x\cM)^\perp$. Since $p(T_x\cM)$ is compatible with the Hodge decomposition \cite{FilipKZ}, so is $p(T_x\cM)^\perp$, because the cup product pairing vanishes on pairs of classes of the same type. Thus, setting
\[
N_x=p(T_x\cM)^\perp\cap H^{1,0}(X),
\]
we have
\[
H^{1,0}(X)=p(T_x\cM)^{1,0}\oplus N_x,\qquad\dim N_x=g-r.
\]

A Beltrami differential on $X$ is a smooth $(0,1)$-form with values in $T_X$, written locally as $\mu(z)\,\frac{d\bar z}{dz}$ in a coordinate $z$. If $|\mu|<1$ everywhere, it determines a complex structure on the surface underlying $X$ whose $(1,0)$-forms are locally the multiples of $dz+\mu(z)\,d\bar z$. By the Dolbeault isomorphism,
\[
H^1(X,T_X)\cong\frac{\{\text{Beltrami differentials on }X\}}{\{\bar\partial v:\ v\text{ a smooth vector field on }X\}},
\]
so every Beltrami differential represents a class in $H^1(X,T_X)$. If $\mu_t$ is a smooth family of Beltrami differentials with $\mu_0=0$, and $X_t$ denotes the surface underlying $X$ with the complex structure determined by $\mu_t$, then under the Kodaira--Spencer identification $T_X\cM_g\cong H^1(X,T_X)$, the tangent vector of $X_t$ at $t=0$ is the class of $\left.\frac{d\mu_t}{dt}\right|_{t=0}$ (see, e.g., \cite[proof of Theorem~1.1]{KahnWright}). We write $\theta=\theta^{1,0}+\theta^{0,1}$ for the decomposition of a smooth complex-valued one form $\theta$ by type.

\begin{lemma}\label{lem:dpi}
Let $\xi\in H^1(X,\Sigma;\bC)$ be represented by $\theta$ as above. Then $d\pi_x(\xi)$ is represented by the Beltrami differential $\theta^{0,1}/\omega$.
\end{lemma}

\begin{proof}
We follow Kahn and Wright \cite[Lemma~2.2 and proof of Theorem~1.1]{KahnWright}. In a local coordinate $z$ on $X\setminus\Sigma$, write $\omega=f\,dz$ and $\theta=A\,dz+B\,d\bar z$. For $t\in\bC$, consider the Beltrami differential
\[
\mu_t=\frac{tB}{f+tA}\,\frac{d\bar z}{dz}.
\]
Since $\theta$ is supported in a compact subset of $X\setminus\Sigma$, on which $\omega$ does not vanish, the function $\theta^{1,0}/\omega$ and the Beltrami differential $\theta^{0,1}/\omega$ are bounded. As
\[
\mu_t=\frac{t\,\theta^{0,1}/\omega}{1+t\,\theta^{1,0}/\omega},
\]
for sufficiently small $t$ the denominator does not vanish, and $\mu_t$ extends smoothly by zero across $\Sigma$ and satisfies $|\mu_t|<1$ everywhere.

Let $X_t$ denote $X$ equipped with the complex structure determined by $\mu_t$. Since
\[
\omega+t\theta=(f+tA)\left(dz+\frac{tB}{f+tA}\,d\bar z\right),
\]
the form $\omega+t\theta$ is of type $(1,0)$ on $X_t$, and it is holomorphic because it is closed. It has no zeros on the support of $\theta$ and agrees with $\omega$ elsewhere, so $(X_t,\omega+t\theta)$ lies in the same stratum. Its relative periods are $\widetilde\omega+t\xi$, so its tangent vector at $t=0$ is $\xi$. Thus $d\pi_x(\xi)$ is represented by
\[
\left.\frac{d\mu_t}{dt}\right|_{t=0}
=\frac{B}{f}\,\frac{d\bar z}{dz}
=\frac{\theta^{0,1}}{\omega}.
\]
\end{proof}

\begin{lemma}\label{lem:annihilate}
For $\alpha\in p(T_x\cM)^{1,0}$ and $\eta\in N_x$, one has $\langle\xi,[\alpha\eta/\omega]\rangle=0$ for every $\xi\in T_x\cM$.
\end{lemma}

\begin{proof}
Fix $\xi\in T_x\cM$, represented by $\theta$ as above. By Lemma~\ref{lem:dpi}, the Kodaira--Spencer class $d\pi_x(\xi)$ is represented by $\mu=\theta^{0,1}/\omega$. Choose a holomorphic curve $x_t$ in $\cM$ with $x_0=x$ and tangent vector $\xi$ at $t=0$. Since $p(T\cM)$ is flat and compatible with the Hodge decomposition, $p(T\cM)^{1,0}$ is a holomorphic subbundle, so we may extend $\alpha$ to a holomorphic section $\alpha_t\in p(T_{x_t}\cM)^{1,0}$ along this curve. By flatness, the Gauss--Manin derivative $\nabla_{\partial_t}\alpha_t$ at $t=0$ lies in $p(T_x\cM)$. 

Although each $\alpha_t$ is holomorphic, this derivative may have a nonzero $(0,1)$-component, since the Hodge decomposition varies along the curve. Because $p(T_x\cM)$ is compatible with the Hodge decomposition, this component lies in $p(T_x\cM)^{0,1}$.

We compute this component using Griffiths' description of the infinitesimal variation of Hodge structure \cite[Theorem~2.3]{Voi13}. It involves the cup product in sheaf cohomology
\[
H^1(X,T_X)\otimes H^0(X,K_X)\xrightarrow{\ \smile\ }H^1(X,T_X\otimes K_X)\cong H^1(X,\mathcal O_X)\cong H^{0,1}(X),
\]
where the first isomorphism is induced by the contraction $T_X\otimes K_X\to\mathcal O_X$. Griffiths' formula states that
\[
\operatorname{pr}^{0,1}\Bigl(\left.\nabla_{\partial_t}\alpha_t\right|_{t=0}\Bigr)=d\pi_x(\xi)\smile\alpha,
\]
where $\operatorname{pr}^{0,1}$ denotes the projection to $H^{0,1}(X)$ and the right side is viewed in $H^{0,1}(X)$ via the isomorphisms above. In the Dolbeault description, where $d\pi_x(\xi)$ is represented by $\mu$, the cup product $d\pi_x(\xi)\smile\alpha$ is represented by $\mu\otimes\alpha$, and the subsequent contraction gives $\mu\lrcorner\alpha$. In a local coordinate on $X\setminus\Sigma$ with $\omega=f\,dz$, $\alpha=a\,dz$ and $\theta^{0,1}=B\,d\bar z$, we have
\[
\mu\lrcorner\alpha=\frac{B}{f}\,\alpha\Bigl(\frac{\partial}{\partial z}\Bigr)\,d\bar z=\frac{aB}{f}\,d\bar z=\frac{\alpha\,\theta^{0,1}}{\omega}.
\]
Consequently,
\[
\operatorname{pr}^{0,1}\Bigl(\left.\nabla_{\partial_t}\alpha_t\right|_{t=0}\Bigr)=\left[\frac{\alpha\,\theta^{0,1}}{\omega}\right]\in p(T_x\cM)^{0,1}.
\]

Finally, since $\alpha\eta/\omega$ has type $(1,0)$, only $\theta^{0,1}$ contributes to the pairing, and since $\theta$ vanishes near $\Sigma$, the form $\alpha\,\theta^{0,1}/\omega$ extends by zero to $X$. Thus
\[
\left\langle\xi,\left[\frac{\alpha\eta}{\omega}\right]\right\rangle=\frac{1}{2\pi i}\int_{X\setminus\Sigma}\frac{\alpha\eta}{\omega}\wedge\theta^{0,1}=\frac{1}{2\pi i}\int_X\eta\wedge\frac{\alpha\,\theta^{0,1}}{\omega}.
\]
The form $\alpha\theta^{0,1}/\omega$ need not be closed. However, $\eta$ is holomorphic and hence closed, so for every smooth function $g$,
\[
d(g\eta)=\partial g\wedge\eta+\bar\partial g\wedge\eta=-\eta\wedge\bar\partial g,
\]
since $\partial g\wedge\eta$ is a $(2,0)$-form and vanishes on $X$. By Stokes' theorem, the integral $\int_X\eta\wedge\frac{\alpha\theta^{0,1}}{\omega}$ depends only on the Dolbeault cohomology class of $\alpha\theta^{0,1}/\omega$. We may therefore replace this form by its harmonic representative, which is closed, and identify the integral with the cup product pairing of $\eta$ with this class. This class lies in $p(T_x\cM)^{0,1}$, and $\eta\in N_x$ is orthogonal to $p(T_x\cM)$, so the integral vanishes.
\end{proof}

\begin{proposition}\label{prop:Vxupper}
If $r<g$, then $\dim V_x\le n-s$.
\end{proposition}

\begin{proof}
Let
\[
P=\operatorname{span}_{\bC}\left\{\frac{\alpha\eta}{\omega}:\alpha\in V_x,\ \eta\in N_x\right\}.
\]
Since $\ord_{z_i}\alpha\ge k_i-1$ and $\ord_{z_i}\omega=k_i$, every element of $P$ has at most simple poles which are all in $\Sigma$.

For $\nu\in P$, the residue theorem gives $\sum_i\Res_{z_i}\nu=0$. Moreover, let $\xi\in W_x$, and let $[c]\in\bC^\Sigma/\bC\mathbf1$ correspond to $\xi$ under the exact sequence above. By Lemma~\ref{lem:annihilate} and linearity we have $\langle\xi,[\nu]\rangle=0$, so Lemma~\ref{lem:pairing} gives
\[
\sum_i c(z_i)\Res_{z_i}\nu=0.
\]
Viewing $W_x$ as a subspace of $\bC^\Sigma/\bC\mathbf1$, it follows that the residue vector of $\nu$ lies in
\[
R=\left\{(\rho_i)\in\bC^n:\sum_i\rho_i=0,\ \sum_i c(z_i)\rho_i=0\text{ for all }[c]\in W_x\right\}.
\]
The pairing $([c],\rho)\mapsto\sum_i c(z_i)\rho_i$ is well defined on vectors $\rho$ whose entries sum to zero. It is nondegenerate, since such a vector is determined by its pairings with the indicator functions of the points $z_i$, and both spaces have dimension $n-1$, so it identifies the hyperplane of such vectors with the dual of $\bC^\Sigma/\bC\mathbf1$. Under this identification, $R$ is the annihilator of $W_x$, so $\dim R=n-1-s$.

We claim that the kernel of the residue map $P\to R$ is $N_x$. As elements of $P$ have at most simple poles, this kernel consists of the holomorphic elements of $P$. Every $\eta\in N_x$ is holomorphic and lies in $P$, since $\eta=\omega\eta/\omega$ and $\omega\in V_x$. Conversely, if $\nu\in P$ is holomorphic, then it is a closed form on $X$, so by Lemma~\ref{lem:pairing}, $\langle\xi,[\nu]\rangle$ is $(2\pi i)^{-1}$ times the cup product pairing of $[\nu]$ with $p(\xi)$ for every $\xi\in T_x\cM$. These pairings vanish by Lemma~\ref{lem:annihilate} and linearity, so $[\nu]$ is orthogonal to $p(T_x\cM)$, and hence it lies in $N_x$. Thus the residue map induces an injection $P/N_x\hookrightarrow R$, and
\[
\dim P\le\dim N_x+n-1-s.
\]

For the reverse bound, let $\ell=\dim V_x$ and $m=\dim N_x$, both positive since $\omega\in V_x$ and $r<g$, and fix $q\in X\setminus\Sigma$. By Gaussian elimination on the leading Taylor coefficients at $q$, we may choose bases $\alpha_1,\ldots,\alpha_\ell$ of $V_x$ and $\eta_1,\ldots,\eta_m$ of $N_x$ with strictly increasing orders of vanishing at $q$. Since $\omega$ does not vanish at $q$, the $\ell+m-1$ forms
\[
\frac{\alpha_1\eta_1}{\omega},\ldots,\frac{\alpha_1\eta_m}{\omega},\frac{\alpha_2\eta_m}{\omega},\ldots,\frac{\alpha_\ell\eta_m}{\omega}
\]
also have strictly increasing orders of vanishing at $q$, so they are linearly independent and $\dim P\ge\ell+m-1$. Combined with the upper bound $\dim P\le m+n-1-s$, this gives $\ell\le n-s$.
\end{proof}

\begin{proof}[Proof of Theorem~\ref{thm:upper}]
By Lemma~\ref{lem:kerdpi} and Proposition~\ref{prop:Vxupper},
\[
\dim\ker d(\pi|_\cM)_x=\dim V_x\le n-s.
\]
To bound the image, consider the span of the holomorphic quadratic differentials $\alpha\eta$ with $\alpha\in p(T_x\cM)^{1,0}$ and $\eta\in N_x$. These two spaces of one forms have dimensions $r\ge1$ and $g-r\ge1$, since $\omega\in p(T_x\cM)^{1,0}$ and $r<g$. Applying the vanishing-order argument from the proof of Proposition~\ref{prop:Vxupper} to the products $\alpha\eta$, we see that this span has dimension at least $r+(g-r)-1=g-1$.

By Serre duality, the pairing between $H^1(X,T_X)$ and $H^0(X,K_X^{\otimes2})$ is perfect \cite[Section~III.7]{Hartshorne}, and up to a nonzero constant, the pairing of the class of a Beltrami differential with a quadratic differential is the integral over $X$ of their product \cite[Section~2]{KahnWright}. Let $\xi\in T_x\cM$ be represented by $\theta$ as above. By Lemma~\ref{lem:dpi}, $d\pi_x(\xi)$ is represented by $\theta^{0,1}/\omega$, so its pairing with $\alpha\eta$ is a nonzero constant multiple of
\[
\int_X\eta\wedge\frac{\alpha\,\theta^{0,1}}{\omega},
\]
which vanishes by Lemma~\ref{lem:annihilate}. Thus these quadratic differentials annihilate the image of $d(\pi|_\cM)_x$. Since $h^0(X,K_X^{\otimes2})=3g-3$, the nondegeneracy of the pairing gives
\[
\dim\operatorname{im}d(\pi|_\cM)_x\le(3g-3)-(g-1)=2g-2.
\]

By the rank--nullity theorem, the bounds for the kernel and the image give
\[
2r+s=\dim T_x\cM\le(n-s)+(2g-2),
\]
so $r+s\le g-1+n/2$. Since $r+s$ is an integer, it follows that $r+s\le g-1+\lfloor n/2\rfloor$.
\end{proof}

\begin{proof}[Proof of Corollary~\ref{cor:fullrel}]
If $s=n-1$ and $r<g$, Theorem~\ref{thm:upper} gives $r\le g-\lceil n/2\rceil$. Now let $\cM$ be a proper invariant subvariety of maximal rel in a principal stratum, where $n=2g-2$ and $s=2g-3$. If $r=g$, then $\cM$ has the same dimension as the stratum and hence equals it by connectedness \cite{KZ03}, contradicting properness. Thus $r<g$, and the bound above gives $r=1$.

Suppose further that $g\ge3$, we show that $\cM$ is defined over $\mathbb Q$. Otherwise, choose a nonidentity embedding $\sigma$ of its field of definition. By Wright's decomposition \cite[Theorem~1.5]{Wri14} and Filip's results \cite{FilipKZ}, the Galois conjugate $p(T\cM)^\sigma$ is a flat symplectic subbundle compatible with the Hodge decomposition, and it intersects $p(T\cM)$ trivially. At $x=(X,\omega)$, set
\[
N'_x=(p(T_x\cM)^\sigma)^\perp\cap H^{1,0}(X).
\]
Since $p(T_x\cM)^\sigma$ has complex dimension two, its $(1,0)$-part is one-dimensional, and hence $\dim N'_x=g-1$. Pick a nonzero $\alpha\in(p(T_x\cM)^\sigma)^{1,0}$ which is not a multiple of $\omega$.

The argument in Lemma~\ref{lem:annihilate} applies more generally to any flat subbundle of $H^1$ over $\cM$ that is compatible with the Hodge decomposition. Applying it to $p(T\cM)^\sigma$, we obtain
\[
\langle\xi,[\alpha\eta/\omega]\rangle=0\qquad\text{for all }\xi\in T_x\cM, \ \eta\in N'_x.
\]
Since $\cM$ has maximal rel, every class in $\ker p$ lies in $T_x\cM$. Applying Lemma~\ref{lem:pairing} to the class of the function equal to $1$ at $z_i$ and $0$ elsewhere, we obtain $\Res_{z_i}(\alpha\eta/\omega)=0$ for every $i$. As $\omega$ has only simple zeros, this form has at most simple poles, and since its residues vanish, it is holomorphic.

The function $\alpha/\omega$ is nonconstant and has only simple poles. It has at least two distinct poles $z_i$ and $z_j$, since a nonconstant meromorphic function with a single simple pole would give an isomorphism $X\cong\mathbb P^1$. As $\alpha\eta/\omega$ is holomorphic for every $\eta\in N'_x$, each such $\eta$ vanishes at $z_i$ and $z_j$, so
\[
N'_x\subset H^0(X,K_X(-z_i-z_j)).
\]
Since $\dim N'_x=g-1$, the Riemann--Roch theorem gives
\[
h^0(X,\mathcal O_X(z_i+z_j))=h^0(X,K_X(-z_i-z_j))+3-g\ge2.
\]
Thus $H^0(X,\mathcal O_X(z_i+z_j))$ contains a nonconstant function, which has at most simple poles at $z_i$ and $z_j$ and no other poles. Since $X$ is not isomorphic to $\mathbb P^1$, this function has degree two, so $X$ is hyperelliptic.

Since $x$ was arbitrary, $\cM$ is contained in the hyperelliptic locus of the principal stratum. Its tangent spaces are then contained in those of this locus, which has rel $g-1$ by Remark~\ref{rem:sharp}, so $s\le g-1$, contradicting $s=2g-3>g-1$. Thus $\cM$ is an arithmetic rank one invariant subvariety, and hence a locus of torus covers.
\end{proof}

Equality in Theorem~\ref{thm:main} holds for at least one invariant subvariety in each stratum whose zeros have orders at most two, and in the principal stratum of genus $g$ it holds for invariant subvarieties of every rank $1\le r\le\lfloor(g+1)/2\rfloor$. For every $g\ge2$ and $1\le n\le2g-2$, equality in Theorem~\ref{thm:upper} holds for some invariant subvariety in some genus $g$ stratum with $n$ zeros, and in the principal stratum it holds in ranks one and $g-1$. We describe these examples in the following remark.

\begin{remark}\label{rem:sharp}
We write $\widetilde{\cQ}(\kappa)$ for the locus of holonomy double covers of surfaces in a stratum $\cQ(\kappa)$ of quadratic differentials. For integers $a,b\ge0$ with $a+b\ge1$, the stratum $\cQ(1^a,2^b,-1^{a+2b+4})$ of quadratic differentials on $\mathbb P^1$ is nonempty \cite{MS93}, and a component of
\[
\widetilde{\cQ}(1^a,2^b,-1^{a+2b+4})\subset\cH(2^a,1^{2b})
\]
has rank $a+b+1$ and rel $b$ by \cite[Lemma~4.2]{AW21}. Since $\sum_i(k_i-1)=a$ in $\cH(2^a,1^{2b})$, equality holds in Theorem~\ref{thm:main}. In particular, $\widetilde{\cQ}(2^b,-1^{2b+4})$ is the hyperelliptic locus in the principal stratum of genus $b+1$.

In the principal stratum of genus $g$, a square-tiled Teichm\"uller curve gives equality in Theorem~\ref{thm:main} in rank one. For $2\le r\le\lfloor(g+1)/2\rfloor$, write
\[
g-1=d(r-1)+k,\qquad d\ge2,\quad 0\le k<r-1.
\]
Let $\mathcal N$ be a component of the locus of connected degree $d$ translation covers of surfaces in $\widetilde{\cQ}(2^{r-1},-1^{2r+2})$ that are simply branched over $2k$ distinct Weierstrass points and unramified elsewhere. To see that such covers exist, take standard generators $a_1,b_1,\ldots,a_r,b_r$ of the fundamental group of the base surface and loops $c_1,\ldots,c_{2k}$ around the branch points, so that $\prod_i[a_i,b_i]\prod_jc_j=1$. Sending $a_1$ to a $d$-cycle, each $c_j$ to the transposition $(1\,2)$, and the remaining generators to the identity respects this relation and gives a transitive action, hence a connected cover with the required branching. By Riemann--Hurwitz, the covers have genus $g$. The base differentials have only simple zeros and do not vanish at the Weierstrass points, so their pullbacks have only simple zeros.

The base locus has rank $r$ and rel $r-1$, and since the Weierstrass points are periodic \cite[p.~3]{AW21}, marking them does not change rank or rel. The tangent space of the covering locus is the pullback of that of the marked base locus \cite[Examples~2.2--2.3]{Wright15}. This pullback is injective on absolute cohomology, and the exact sequences of the pairs show that it is also injective on relative cohomology, since every zero and every branch point downstairs has a zero among its preimages. Hence $\mathcal N$ has rank $r$ and rel $r-1$, giving equality in Theorem~\ref{thm:main}.

For Theorem~\ref{thm:upper} in the principal stratum of genus $g$, take a component of the locus of connected double covers of flat tori branched over $2g-2$ distinct points. This locus has rank one and dimension $2g-1$, hence rel $2g-3$, so equality holds. An example of rank $g-1$ is given by the case $n=2g-2$ below.

Finally, fix $g\ge2$ and $1\le n\le2g-2$. For $(g,n)=(2,1)$, a square-tiled Teichm\"uller curve in $\cH(2)$ gives equality in Theorem~\ref{thm:upper}. Otherwise, write $n=2p$ or $n=2p+1$. In the first case, choose positive integers $a_1,\ldots,a_p$ with $\sum_ja_j=g-1$, and in the second case choose positive integers $b,a_1,\ldots,a_p$ with $b+\sum_ja_j=g-1$. Consider a component of
\[
\widetilde{\cQ}(2a_1,\ldots,2a_p,-1^{2g-2})\quad\text{or}\quad\widetilde{\cQ}(2b-1,2a_1,\ldots,2a_p,-1^{2g-3})
\]
according to the parity of $n$. These are strata of quadratic differentials on tori, and they are nonempty \cite{MS93}. The double covers have genus $g$ and $n$ zeros, also these loci have rank $g-1$ and rel $p=\lfloor n/2\rfloor$, so equality holds in Theorem~\ref{thm:upper}. For $n=2g-2$, taking all $a_j=1$ gives the locus $\widetilde{\cQ}(2^{g-1},-1^{2g-2})$ in the principal stratum, which has rank and rel both equal to $g-1$.
\end{remark}


\begin{thebibliography}{BCGGM19}

\bibitem[Api18]{Api18}
P.~Apisa, \emph{$\mathrm{GL}_2\mathbb{R}$ orbit closures in hyperelliptic components of strata},
Duke Math. J. \textbf{167} (2018), no.~4, 679--742.

\bibitem[AW21]{AW21}
P.~Apisa and A.~Wright, \emph{Marked points on translation surfaces},
Geom. Topol. \textbf{25} (2021), no.~6, 2913--2961.

\bibitem[AW23]{AW23}
P.~Apisa and A.~Wright, \emph{High rank invariant subvarieties},
Ann. of Math. (2) \textbf{198} (2023), no.~2, 657--726.

\bibitem[AN16]{AN16}
D.~Aulicino and D.-M.~Nguyen, \emph{Rank two affine submanifolds in $\mathcal{H}(2,2)$ and $\mathcal{H}(3,1)$},
Geom. Topol. \textbf{20} (2016), no.~5, 2837--2904.

\bibitem[AN20]{AN20}
D.~Aulicino and D.-M.~Nguyen, \emph{Rank 2 affine manifolds in genus~3},
J. Differential Geom. \textbf{116} (2020), no.~2, 205--280.

\bibitem[ANW16]{ANW16}
D.~Aulicino, D.-M.~Nguyen and A.~Wright, \emph{Classification of higher rank orbit closures in $\mathcal{H}^{\mathrm{odd}}(4)$},
J. Eur. Math. Soc. (JEMS) \textbf{18} (2016), no.~8, 1855--1872.

\bibitem[AEM17]{AEM17}
A.~Avila, A.~Eskin and M.~M\"oller, \emph{Symplectic and isometric $\mathrm{SL}(2,\mathbb{R})$-invariant subbundles of the Hodge bundle},
J. Reine Angew. Math. \textbf{732} (2017), 1--20.

\bibitem[BCGGM19]{BCGGM19}
M.~Bainbridge, D.~Chen, Q.~Gendron, S.~Grushevsky and M.~M\"oller, \emph{Strata of $k$-differentials},
Algebr. Geom. \textbf{6} (2019), no.~2, 196--233.

\bibitem[CDF23]{CDF23}
G.~Calsamiglia, B.~Deroin and S.~Francaviglia, \emph{A transfer principle: from periods to isoperiodic foliations},
Geom. Funct. Anal. \textbf{33} (2023), no.~1, 57--169.

\bibitem[CW24]{CW24}
J.~Chaika and B.~Weiss, \emph{On the ergodic theory of the real Rel foliation},
Forum Math. Pi \textbf{12} (2024), e7.

\bibitem[EFW18]{EFW18}
A.~Eskin, S.~Filip and A.~Wright, \emph{The algebraic hull of the Kontsevich--Zorich cocycle},
Ann. of Math. (2) \textbf{188} (2018), no.~1, 281--313.

\bibitem[EM18]{EM18}
A.~Eskin and M.~Mirzakhani, \emph{Invariant and stationary measures for the $\mathrm{SL}(2,\mathbb{R})$ action on moduli space},
Publ. Math. Inst. Hautes \'Etudes Sci. \textbf{127} (2018), 95--324.

\bibitem[EMM15]{EMM15}
A.~Eskin, M.~Mirzakhani and A.~Mohammadi, \emph{Isolation, equidistribution, and orbit closures for the $\mathrm{SL}(2,\mathbb{R})$ action on moduli space},
Ann. of Math. (2) \textbf{182} (2015), no.~2, 673--721.

\bibitem[Fil16a]{FilipKZ}
S.~Filip, \emph{Semisimplicity and rigidity of the Kontsevich--Zorich cocycle},
Invent. Math. \textbf{205} (2016), no.~3, 617--670.

\bibitem[Fil16b]{FilipMHS}
S.~Filip, \emph{Splitting mixed Hodge structures over affine invariant manifolds},
Ann. of Math. (2) \textbf{183} (2016), no.~2, 681--713.

\bibitem[Har77]{Hartshorne}
R.~Hartshorne, \emph{Algebraic Geometry},
Graduate Texts in Mathematics \textbf{52}, Springer, New York, 1977.

\bibitem[KW22]{KahnWright}
J.~Kahn and A.~Wright, \emph{Hodge and Teichm\"uller},
J. Mod. Dyn. \textbf{18} (2022), 149--160.

\bibitem[KZ03]{KZ03}
M.~Kontsevich and A.~Zorich, \emph{Connected components of the moduli spaces of Abelian differentials with prescribed singularities},
Invent. Math. \textbf{153} (2003), no.~3, 631--678.

\bibitem[Lee24]{Lee24}
J.~M.~Lee, \emph{Introduction to Complex Manifolds},
Graduate Studies in Mathematics 244, American Mathematical Society,
Providence, RI, 2024.

\bibitem[MS93]{MS93}
H.~Masur and J.~Smillie, \emph{Quadratic differentials with prescribed singularities and pseudo-Anosov diffeomorphisms},
Comment. Math. Helv. \textbf{68} (1993), no.~2, 289--307.

\bibitem[McM07]{McM07}
C.~T.~McMullen, \emph{Dynamics of $\mathrm{SL}_2(\mathbb{R})$ over moduli space in genus two},
Ann. of Math. (2) \textbf{165} (2007), no.~2, 397--456.

\bibitem[MW18]{MW18}
M.~Mirzakhani and A.~Wright, \emph{Full-rank affine invariant submanifolds},
Duke Math. J. \textbf{167} (2018), no.~1, 1--40.

\bibitem[M\"ol08]{MollerLinear}
M.~M\"oller, \emph{Linear manifolds in the moduli space of one-forms},
Duke Math. J. \textbf{144} (2008), no.~3, 447--487.

\bibitem[Voi13]{Voi13}
C.~Voisin, \emph{Hodge loci}, in \emph{Handbook of Moduli. Vol.~III},
Adv. Lect. Math. (ALM) \textbf{26}, Int. Press, Somerville, MA, 2013, 507--546.

\bibitem[Win21]{Win21}
K.~Winsor, \emph{Dynamics of the absolute period foliation of a stratum of holomorphic 1-forms},
arXiv:2109.12669.

\bibitem[Win24]{Win24}
K.~Winsor, \emph{Saturated orbit closures in the Hodge bundle},
Int. Math. Res. Not. IMRN \textbf{2024}, no.~6, 4710--4737.

\bibitem[Win26]{Win26}
K.~Winsor, \emph{Isoperiodic forms and invariant subvarieties of connected strata},
arXiv:2606.01587.

\bibitem[Wri14]{Wri14}
A.~Wright, \emph{The field of definition of affine invariant submanifolds of the moduli space of abelian differentials},
Geom. Topol. \textbf{18} (2014), no.~3, 1323--1341.

\bibitem[Wri15]{Wright15}
A.~Wright, \emph{Translation surfaces and their orbit closures: an introduction for a broad audience},
EMS Surv. Math. Sci. \textbf{2} (2015), no.~1, 63--108.

\end{thebibliography}
\end{document}